# Modern efficient numerical approaches to regularized regression problems in application to traffic demands matrix calculation from link loads

*Anikin Anton*

(Institute of System Dynamics and Control Theory, Siberian Branch of the Russian Academy of Sciences, Irkutsk) anton.anikin@htower.ru

*Dvurechensky Pavel*

(IITP RAS, WIAS Berlin, PreMoLab MIPT) pavel.dvurechensky@gmail.com

*Gasnikov Alexander*

(IITP RAS, PreMoLab MIPT, HSE) avgasnikov@gmail.com

*Golov Andrey*

(PreMoLab MIPT) and_go@mail.ru

*Gornov Alexander*

(Institute of System Dynamics and Control Theory, Siberian Branch of the Russian Academy of Sciences, Irkutsk) gornov.a.yu@gmail.com

*Maximov Yuri*

(IITP RAS, HSE, PreMoLab MIPT) yuramaksimov@gmail.com

*Mendel Mikhail*

(PreMoLab MIPT) mendelm@mail.ru

*Spokoiny Vladimir*

(WIAS Berlin, IITP RAS, PreMoLab MIPT) spokoiny@wias-berlin.de

**Abstract**

In this work we collect and compare to each other many different numerical methods for regularized regression problem and for the problem of projection on a hyperplane. Such problems arise, for example, as a subproblem of demand matrix estimation in IP-networks. In this special case matrix of affine constraints has special structure: all elements are 0 or 1 and this matrix is sparse enough. The number of rows can be $10^5$ and the number of columns can be $10^8$. So we have to deal with huge-scale convex optimization problem of special type. Using the properties of the problem we try "to look inside the black-box" and to see how the best modern methods work when applied to this problem.

**Introduction.**

In the problem of traffic demand matrix estimation the goal is to recover traffic demand matrix represented as a vector $x \geq 0$ from known route matrix $A$ (the element $[A]_{i,j}$ is equal 1 iff the demand with number $j$ goes through link with number $i$ and equals 0 otherwise) and link loads $b$ (amount of traffic which goes through every link). This leads to the problem of finding the solution of linear system $Ax = b$. Also we assume that we have some $x_g \geq 0$ which reflects our prior assumption about $x$. Thus we consider $x$ to be a projection of $x_g$ on a simplex-type set $\{x \geq 0 : Ax = b\}$

$$\min_{\substack{Ax=b \\ x \geq 0}} \left\{ g(x) := \|x - x_g\|_2^2 \right\} = \min_{\substack{\|Ax-b\|_2^2 \leq 0 \\ x \geq 0}} g(x).$$

Slater's relaxation of this problem leads to the problem (denote $x_*$ the solution of this problem)

$$\|x - x_g\|_2^2 \to \min_{\substack{\|Ax-b\|_2^2 \leq \varepsilon^2 \\ x \geq 0}}.$$

This problem can be reduced to the problem (unfortunately without explicit dependence $\bar{\lambda}(\varepsilon)$)

$$\tilde{f}(x) = \|x - x_g\|_2^2 + \bar{\lambda} \|Ax - b\|_2^2 \to \min_{x \geq 0},$$

where $\bar{\lambda}$ – dual multiplier to the convex inequality $\|Ax - b\|_2^2 \leq \varepsilon^2$. One might expect that $\bar{\lambda} \gg \|x_* - x_g\|_2^2 / \varepsilon^2$, but in reality $\bar{\lambda}$ can be chosen much smaller ($\bar{\lambda} \sim \varepsilon^{-1} - \varepsilon^{-2}$, see Part 2) if we restrict ourselves only by approximate solution. Let's reformulate the problem

$$f(x) = \|Ax - b\|_2^2 + \lambda \|x - x_g\|_2^2 \to \min_{x \geq 0},$$

where $\lambda = \bar{\lambda}^{-1}$. The last two problem statements can be considered as problems of Bayesian parameter estimation [1]: One measure the vector $b$ with some random error $\xi \in N(0, \sigma^2 I)$ and tries to find such vector $x \geq 0$ that satisfies the linear system $Ax = b$ assuming that this vector is also random with prior distribution $x \in N(x_g, \tilde{\sigma}^2 I)$. So the model of the data is the following.

$$b = Ax + \xi, \ \xi \in N(0, \sigma^2 I), \text{ prior on } x \in N(x_g, \tilde{\sigma}^2 I).$$

So the functions $\tilde{f}(x)$ and $f(x)$ introduced above now within multipliers are minus log likelihood for this model (with $\bar{\lambda} = \tilde{\sigma}^2 / \sigma^2$, $\lambda = \sigma^2 / \tilde{\sigma}^2$). Hence one can consider e.g. the second minimization problem as a Bayesian estimation problem [1] or as a Penalized Maximum Likelihood Estimation (V. Spokoiny, 2012 [2]).

In this paper we consider not only Euclidian projection. The second natural choice is Kullback–Leibler "projection". We also consider a problem of finding a sparse solution of the system $Ax = b$ which leads to LASSO-type problem.

The main result of the work is overview of modern approaches for the numerical solution of the mentioned above problems. The main practical motivation for us IP-traffic analysis (see e.g. [3], [4]). We also slightly generalize some known results.

- We describe fast gradient method for composite optimization problems in which entropy function is considered as composite term;
- We propose new variant of gradient-free Powell's type method based on calculation of function values at three points at each iteration;
- We describe random coordinate descent for our particular case ($A$ has all elements equal 0 or 1);
- We propose new estimates for dual accelerated random coordinate ascent that allows to use sparsity of the problem;
- We propose a technique that allows to reduce the projection problem to the small number of regularized regression type problem.

**Part 1. Regularized Regression Approach**

Consider the following problem (instead of $\varphi(x)$ we can considered many other sum-type function of scalar product of rows some sparse matrix and $x$ – most of the results below can be generalized in this direction)

$$f(x) = \underbrace{\frac{1}{2}\|Ax - b\|_2^2}_{\varphi(x)} + g(x) \to \min_{x \in Q}, \qquad (1)$$

where $A$ – is a matrix of size $m \times n$ with elements equal 0 or 1. We assume that the matrix $A$ is sparse with number of non-zero elements $nnz(A)$. Set $s = nnz(A)/n \ll m$, $\tilde{s} = sn/m$. By the solution of the problem (1) we will mean such vector (generally speaking a random vector if the method is randomized) $x^N$, that[1]

---

[1] Since all the situation we considered below can be treated in a strongly convex environment we assume that the high probability deviations bounds can be obtained from the Markov inequality

$$P\left(f\left(x^{N(\varepsilon^2\sigma)}\right) - f_* \geq \varepsilon^2\sigma/\sigma\right) \leq \frac{E\left[f\left(x^{N(\varepsilon^2\sigma)}\right)\right] - f_*}{\varepsilon^2\sigma/\sigma} \leq \sigma, \; N(\varepsilon) \sim \sqrt{L/\lambda}\ln\left(LR^2/\varepsilon\right).$$

$$E\left[f\left(x^N\right)\right] - f_* \leq \varepsilon^2. \tag{2}$$

Where the expectation is taken with respect to all randomness in the method.

Possible cases for choice of $g(x)$ are:

1. (Ridge Regression [1] / Tomogravity model [3])

$$g(x) = \lambda \|x - x^g\|_2^2, \ Q = \mathbb{R}_+^n;$$

2. (Mimimal mutual information model [4])[2]

$$g(x) = \lambda \sum_{k=1}^n x_k \ln\left(x_k / x_k^g\right), \ x_k^g \in Q = S_n(R) = \left\{x \geq 0 : \sum_{k=1}^n x_k = R\right\}.$$

All the results below can be generalized for $Q = \breve{S}_n(R) = \left\{x \geq 0 : \sum_{k=1}^n x_k \leq R\right\}$;

3. (LASSO [1])

$$g(x) = \lambda \|x\|_1, \ Q = \mathbb{R}_+^n;$$

Parameter $\lambda$ is a structural parameter. Note that in the third case one often needs to solve the problem (1) for different $\lambda$. For the first and second cases one usually chooses $\lambda \sim \chi \varepsilon^2$ (how to choose $\chi$ will be described in Part 2).

We use the following notations:

$\sigma_{\max}(A)$ – maximal eigenvalue of the matrix $A^T A$, note that

$$\sigma_{\max}(A) = \lambda_{\max}\left(A^T A\right) \leq \operatorname{tr}\left(A^T A\right) = \sum_{k=1}^n \left\|A^{\langle k \rangle}\right\|_2^2 = nnz(A) = sn,$$

where $A^{\langle k \rangle}$ – $k$-th column of matrix $A$;

$\max_{k=1,\ldots,n} \left\|A^{\langle k \rangle}\right\|_2^2 \leq m$;

$R_2^2 = \frac{1}{2}\|x^0 - x_*\|_2^2$, where $x^0$ – starting point, $x_*$ – solution of (1);

$\tilde{O}(\ ) = O(\ )$ up to a logarithmic factor.

In the table below one can find complexity estimates (mathematical expectation of the total number of flops operations needed for finding solution of the problem (1) in sense (2)) for different algorithms applied to the problem (1) with different choice of $g(x)$. We marked by *star the situations which are new in some extent.

---

So we just have to make $\ln(\sigma^{-1})$-times additional iterations to have $\geq 1-\sigma$ probability guarantee.

[2] Below for simplicity we assume that the components of the vector $x^g$ are equal.

| algorithm/ model | Ridge Regression / Tomogravity model | Mimimal mutual information model | LASSO |
|---|---|---|---|
| Conjugate Gradients Method and different modifications [5] | $\min\left\{\begin{array}{l}\tilde{O}\left(sn\sqrt{\dfrac{\sigma_{\max}(A)}{\lambda_1}}\right)\\ O\left(sn\sqrt{\dfrac{\sigma_{\max}(A)R_2^2}{\varepsilon^2}}\right)\end{array}\right\}^3$ | Not applicable | Not applicable |
| Composite FGM [6] | $\min\left\{\begin{array}{l}\tilde{O}\left(sn\sqrt{\dfrac{\sigma_{\max}(A)}{\lambda_1}}\right)\\ O\left(sn\sqrt{\dfrac{\sigma_{\max}(A)R_2^2}{\varepsilon^2}}\right)\end{array}\right\}$ | $\min\left\{\begin{array}{l}\tilde{O}\left(sn\sqrt{\dfrac{\max_{k=1,\ldots,n}\|A^{\langle k\rangle}\|_2^2 R}{\lambda_2}}\right)\\ \tilde{O}\left(sn\sqrt{\dfrac{\max_{k=1,\ldots,n}\|A^{\langle k\rangle}\|_2^2 R^2}{\varepsilon^2}}\right)\end{array}\right\}*$ | $O\left(sn\sqrt{\dfrac{\sigma_{\max}(A)R_2^2}{\varepsilon^2}}\right)$ |
| 3 point plane method of [5] Gornov–Anikin–Powell | $\min\left\{\begin{array}{l}\tilde{O}\left(sn\dfrac{s}{\lambda_1}\right)\\ \tilde{O}\left(sn\dfrac{sR_2^2}{\varepsilon^2}\right)\end{array}\right\}*$ | Not applicable | $\tilde{O}\left(sn\dfrac{sR_2^2}{\varepsilon^2}\right)*$ |
| RCD APPROX / ALPHA [7], [8] | $\min\left\{\begin{array}{l}\tilde{O}\left(sn\sqrt{\dfrac{s}{\lambda_1}}\right)*\\ O\left(sn\sqrt{\dfrac{sR_2^2}{\varepsilon^2}}\right)\end{array}\right\}$ | Not applicable | $O\left(sn\sqrt{\dfrac{sR_2^2}{\varepsilon^2}}\right)$ |
| Dual RCA [9],[10] | $\min\left\{\begin{array}{l}\tilde{O}\left(\tilde{s}m\sqrt{\dfrac{\tilde{s}}{\lambda_1}}\right)\\ \tilde{O}\left(\tilde{s}m\sqrt{\dfrac{\tilde{s}R_2^2}{\varepsilon^2}}\right)\end{array}\right\}*$ | $\min\left\{\begin{array}{l}\tilde{O}\left(mn\sqrt{\dfrac{R}{\lambda_2}}\right)\\ \tilde{O}\left(mn\sqrt{\dfrac{R^2}{\varepsilon^2}}\right)\end{array}\right\}*$ | Not applicable |

It should be mentioned that parameter $R_2^2$ is typically unknown a priory. But this parameter is included in step size policy of many methods in the table above, moreover this parameter is included in the stopping rule criteria of many methods, since we don't know optimal value of the problem $f_*$ (especially in case 3). So what we have to do? The basic idea is proper regularization of the problem (in case 1 the complexity of iteration doesn't change in case 3 it doesn't change too because of the shrinkage operator trick [11]). We choose such $R_2$ that $\varepsilon^2/R_2^2 > 2\lambda_1$ and choose regularization parameter $\mu = \varepsilon^2/(2R_2^2)$. We believe that we guess the true value of $R_2$. We calculate according to this $R_2^2$ prescribed number of iteration $N(R_2,\varepsilon)$ (we have explicitly formulas for that). Since we have strictly convex case we can formulate the stopping rule in terms of gradient

---

[3] Strictly speaking, we know the proof of this estimation (see [5]) only in the case when $Q = \mathbb{R}^n$.

$$f(x^N) - f_* \leq f(x^N) - \min_{x \in \mathbb{R}_+^n} \left\{ f(x^N) + \langle \nabla f(x^N), x - x^N \rangle + \frac{\mu}{2} \|x - x^N\|_2^2 \right\} =$$

$$= \max_{x \in \mathbb{R}_+^n} \left\{ \langle \nabla f(x^N), x^N - x \rangle - \frac{\mu}{2} \|x - x^N\|_2^2 \right\} \leq \varepsilon^2.$$

We verify our guess according to the performance of this stopping rule (due to the linear rate of convergence the roughness of the stopping rule doesn't matter – it will costs for us only additional logarithmic factor in final estimation). If the stopping rule is fulfilled, all right: we stop! Otherwise we put $R_2 := 2R_2$ and restart procedure. We have to do no more than $\log_2(2\lambda_1 R_2^2 / \varepsilon^2)$ such restarts ( $R_2$ in this formula is starting $R_2$ ).

**Remark 1 (FGM for composite problems).** Minimal mutual information model in strongly convex case (case 2) fits well to the framework of composite optimization [6] since one can consider $g(x) = \lambda \sum_{k=1}^n x_k \ln(x_k / x_k^g)$ as the composite term. It is sufficient that $g(x)$ is $\lambda$-strongly convex in 1-norm. Smoothness of $g(x)$ is not required. We use composite FGM method with 1-norm in primal $x$-space and prox-function $d(x) = \ln n + \sum_{k=1}^n x_k \ln x_k$ in not strongly convex case and

$$d(x) = \frac{1}{2(a-1)} \|x\|_a^2 \quad \text{with} \quad a = \frac{2 \log n}{2 \log n - 1}$$

in strongly convex case [6], [12]. In not strongly convex case we can calculate the new point according to explicit formulas (exponential weighting), because prox-term and composite one both are entropy-type in strongly convex case this also can be done effective (see text after Remark 4). Note that in not strongly convex case we have to use $R^2 \ln n$ instead of $R^2$. In strongly convex case we also have to replace $R$ by $R \ln n$ because of the arguments in [13]. Here we have an example when non-Euclidian prox-structure in strongly convex cases gave more benefits than Euclidian one.

**Remark 2 (Powell's type method).** We propose the following generalization of the classic gradient-free method of Powell [5]. Choose a random ort $e_i$. Calculate the value of the function at the points $x^k \pm e_i$, where $x^k$ is a current position. Note that we can find in a fast manner minimum of parabolic function (difficult in calculation because of the size) and simple composite on the line $y = x^k + te_i$, $t \in \mathbb{R}$, because we can determine the 1-parabolic function in unique manner from the values in 3 points ($x^k \pm e_i$ and $x^k$). The only problem is to recalculate the parabolic function in 3 points ($x^k \pm e_i$ and at minimum at the line $x^{k+1}$). Assume that we've already calculated $Ax^k$, since we consider ort $e_i$, we can reduce the problem of recalculation $f(x)$ to the calculation $Ae_i$, that can be done by $O(s)$ flops. So the whole iteration costs $\tilde{O}(s)$. In this situation we have the rate of convergence corresponds to the random primal coordinate descent method [5]. In practice this method works better in many situations than Conjugate Gradients Method.

**Remark 3 (Estimates for FGM-type methods).** Now we try to explain how these estimates were obtained. First two rows of the table have the following form [6], [12]

$$\underbrace{\tilde{O}(sn)}_{\substack{\text{costs of one iteration, the main} \\ \text{part is calculation of full gradient,} \\ \text{that is calculation of } A^T \cdot (Ax)}} \cdot \min\left\{ \underbrace{\tilde{O}\left(\sqrt{L_p/\tilde{\lambda}_p}\right)}_{\substack{\text{the number of iterations,} \\ \text{according to the FGM in} \\ \text{strongly convex case}}}, \underbrace{O\left(\sqrt{L_p R_p^2/\tilde{\varepsilon}}\right)}_{\substack{\text{the number of iterations,} \\ \text{according to the FGM in} \\ \text{non strongly convex case}}} \right\},$$

where $\tilde{\varepsilon}$ – precision we'd like to have. We use $\tilde{\varepsilon} = \varepsilon^2$ (see (2)); $R_p^2$ – Bregman divergence (see [6], [12]) between the starting point and the solution of problem (1) in case when we choose $p$-norm in primal $x$-space (for example in $p = 2$ we have $R_2^2$, introduced above); $\tilde{\lambda}_p$ – constant of strongly convexity $f(x)$ in $p$-norm (in case 1 $p = 2$, $\tilde{\lambda}_2 = \lambda$ and in case 2 $p = 1$, $\tilde{\lambda}_1 = \lambda/R$); $L_p$ – Lipschitz constant of the gradient of $\varphi(x)$:

$$L_p = \max_{x \in Q} \max_{\|h\|_p \leq 1} \left\langle h, \left\|\frac{\partial^2 \varphi}{\partial x_i \partial x_j}\right\| h \right\rangle = \max_{\|h\|_p \leq 1} \left\langle h, A^T A h \right\rangle = \max_{\|h\|_p \leq 1} \left\langle Ah, Ah \right\rangle = \max_{\|h\|_p \leq 1} \|Ah\|_2^2.$$

In cases 1, 3 $p = 2$, $L_2 = \lambda_{\max}\left(A^T A\right) \stackrel{def}{=} \sigma_{\max}(A)$ and in case 2 $p = 1$, $L_1 = \max_{k=1,\ldots,n} \|A^{\langle k \rangle}\|_2^2$.

**Remark 4 (Estimates for Random Coordinate Descent (RCD) methods).** Let's compare FGM-type estimates to its RCD counterparts [7], [8] (see rows 3, 4 of the table)[4]

$$\underbrace{\tilde{O}(s)}_{\substack{\text{costs of one iteration,} \\ \text{the main part is recalculation} \\ \text{of component of gradient}}} \cdot \underbrace{n}_{\substack{\text{payment for} \\ \text{calculation only random} \\ \text{component of gradient}}} \cdot \min\left\{ \underbrace{\tilde{O}\left(\sqrt{\bar{L}_p/\tilde{\lambda}_p}\right)}_{\substack{\text{the number of iterations,} \\ \text{according to the FGM in} \\ \text{strongly convex case}}}, \underbrace{O\left(\sqrt{\bar{L}_p R_p^2/\tilde{\varepsilon}}\right)}_{\substack{\text{the number of iterations,} \\ \text{according to the FGM in} \\ \text{non strongly convex case}}} \right\},$$

where $\bar{L}_p$ is, roughly speaking, the average Lipschitz constant of gradient of $\varphi(x)$:

$$\bar{L}_p^{1/2} = \max_{x \in Q} \frac{1}{n} \sum_{k=1}^n \left\langle e_k, \left\|\frac{\partial^2 \varphi}{\partial x_i \partial x_j}\right\| e_k \right\rangle^{1/2} = \frac{1}{n} \sum_{k=1}^n \|A^{\langle k \rangle}\|_2 \leq \sqrt{s}.$$

Here we considered only the case $p = 2$ and non strongly convex situations [7], [8] with separable composite and set $Q$. But Peter Richtarik announced to us that he has obtained these results also in strongly convex case in August, 2014 (it seems that these results can be obtained by restart technique [13]). Generalization to non Euclidian set up or(and) non separable structure of composite term and set $Q$ to the best of our knowledge hasn't been made until now.

If we assume that $Q$ is formed by a few $r$ affine restrictions (or some others separable convex inequalities), we can insert them with Lagrange multipliers in the goal function. Then we can solve new problem with fixed multipliers with the same complexity. At the same time we can consider the dual problem which has small dimension. To calculate (super-)subgradient of the goal function in the dual problem we have to solve primal problem ($\tilde{\varepsilon}$ solution in terms of function value of the primal problem gives us $\tilde{\varepsilon}$-subgradient of the dual problem [5]). If we use ellipsoids method [14], [15], we can find in a fast manner (since the dimension is small) solution of the dual problem with accuracy $\varepsilon = \mathrm{O}(\tilde{\varepsilon})$ (in terms of dual function value). With appropriate choice of the method for the dual problem (ellipsoids method is proper) one can obtain the solution of the primal problem with the same accuracy in terms of primal function value $\varepsilon$ [15]. This trick (described in 5.3.3.3 [14] and item [15]) can be used in many other contexts (see Part 2).

---

[4] Here we have not trivial assumption about possibility of recalculation of gradient's components (we also have to make a proper step by coordinate-vise version of composite variant of FGM) after $O(s)$ flops. But for many interesting cases (especially problem formulation (1)) it is possible indeed, see p. 16, 17 [7].

So the main advantage of RCD methods consists in change of worth-case Lipschitz constant of gradient in complexity estimates to its average counterpart. This average Lipschitz constant can be much smaller, since (case 1) typically $s \ll \sigma_{\max}(A)$ since

$$\sigma_{\max}(A) = \max_{k=1,\ldots,n} \lambda_k(A^T A), \quad \sum_{k=1}^{n} \lambda_k(A^T A) = \mathrm{tr}(A^T A) = sn.$$

**Remark 5 (Estimates for Dual Random Coordinate Ascend (RCA) methods).** First of all let's form the dual problem. Denote by $A_k$ – $k$-th row of matrix $A$, $\sigma_k(z) = \frac{1}{2}(z - b_k)^2$. Then we have (see Sion–Kakutani minimax theorem at the very end of the book [14])

$$\min_{x \in Q}\left\{\frac{1}{2}\|Ax - b\|_2^2 + g(x)\right\} = \min_{x \in Q}\left\{\sum_{k=1}^{m} \sigma_k(A_k x) + g(x)\right\} =$$

$$= \min_{\substack{x \in Q \\ f = Ax}}\left\{\sum_{k=1}^{m} \sigma_k(f_k) + g(x)\right\} = \min_{\substack{x \in Q \\ f = Ax, f'}} \max_y \left\{\langle f - f', y\rangle + \sum_{k=1}^{m} \sigma_k(f'_k) + g(x)\right\} =$$

$$= \max_{y \in \mathbb{R}^m}\left\{-\max_{\substack{x \in Q \\ f = Ax}}\{\langle f, y\rangle - g(x)\} - \max_{f'}\left\{\langle f', y\rangle - \sum_{k=1}^{m} \sigma_k(f'_k)\right\}\right\} =$$

$$= \max_{y \in \mathbb{R}^m}\left\{-\max_{x \in Q}(\langle -A^T y, x\rangle - g(x)) - \sum_{k=1}^{m} \max_{f'_k}(f'_k y_k - \sigma_k(f'_k))\right\} =$$

$$= \max_{y \in \mathbb{R}^m}\left\{-g^*(-A^T y) - \sum_{k=1}^{m} \sigma_k^*(y_k)\right\} = -\min_{y \in \mathbb{R}^m}\left\{g^*(-A^T y) + \sum_{k=1}^{m} \sigma_k^*(y_k)\right\},$$

where we have explicit expressions for $g^*$, $\sigma_k^*$ and its gradients. Moreover we have explicit dependence of feasible $\bar{x}(y) \in Q$. If $y_*$ is an optimal solution of this dual problem, then $x_* = \bar{x}(y_*)$. Due to duality properties we also have that $\sum_{k=1}^{m} \sigma_k^*(y_k)$ is 1-strongly convex in 2-norm in dual $y$-space and $g^*(-A^T y)$ has Lipschitz constant of gradient in 2-norm equal to $\sigma_{\max}(A^T)/\lambda \leq \tilde{s}m/\lambda = sn/\lambda$ in case 1 and $\max_{k=1,\ldots,n}\|A^{\langle k\rangle}\|_2^2 R/\lambda$ in case 2 (see [9] and theorem 1 in [12]). We can use RCD for the dual problem multiplies by "-1" and use the approach briefly described in the previous remark. It is worth noting that one has a possibility in case 1 to use recalculation at each iteration to obtain complexity of one iteration $\tilde{O}(\tilde{s})$. Unfortunately, in case 2 we can only obtain complexity of one iteration $O(m)$. In this case we also have average Lipschitz constant ($A_{ij} \in \{0,1\} - (i,j)$ element of matrix $A$)

$$\overline{L}^{1/2} \leq 2 \max_{p \in S_n(1)} \frac{1}{m} \sum_{i=1}^{m} \left( \sum_{j=1}^{n} A_{ij}^2 p_j \right)^{1/2} \leq 2.$$

Note that the dual problem is unconstrained. So this is one of the ways to work with not separable constraints in primal problem (but we have payment for that – dual functional isn't still separable, so in sparse case we have lack of possibility to use sparsity for accelerated methods).

It should be mentioned that in the dual approach we also have to restrict ourselves only by Euclidian prox-structure in the dual space and we have to consider strongly convex (concave) case. Except for the construction described above (when we explicitly use strongly convexity of the primal problem[5]) we also have the following reason for that. We'd like to calculate the solution of the primal problem without technique of the work [15], because in this case one iteration will costs $\tilde{O}(n)$ (however there are some exceptions [17]). So the only way to do it is the convergence of the iteration process in the dual space not only in function but also in an argument. Strongly convexity is the only simple way to guarantee it in a friendly computational manner. So we have to consider dual problem as strongly convex (concave) problem.

We are interested in traffic applications [3] in which typically $m \ll n$ (in the machine learning applications the situation is typically inverse $m \gg n$ [1], [18]). So we should use primal RCD.

**Remark 6 (Parallel and distributed computations).** RCD can be fully parallelized in $O(n/\tilde{s})$ processors according to [7], [8]. Typically that method with whole calculation of the gradient can be parallelized in computation of $Ax$ or $A^T y$. Our problem formulation is also well suited for distributed calculation (see ADMM by S. Boyd et al. [18] and solver based on it; see also works of P. Richtarik [19]).

**Remark 7 (Accuracy).** In cases 1, 2 we expect to have such a situation for projection problem when estimates in strongly convex case seems to be close enough to estimates in not strongly convex (arguments at each min in the table above are close to each other). That is in some situations it doesn't matter to use strongly convexity or not – the rates of convergence up to a logarithmic factor are the same. But even in this situation there is a difference. The difference is the following: in the strongly convex case we have guarantee for convergence in argument, but in not strongly convex case we are able to use widely variety of prox-structures. One should also say that in real application we have to choose $\varepsilon$ according to initial discrepancy. So we have to work with relative precision. This fact allows us to fix some level of the relative precision (we choose 0.01, i.e. 1%) and tie the stopping rule of the method to the performance of this criterion.

---

[5] If it isn't a true we can use proper regularization of the primal problem (coincides with composite term up to a multiplicative constant) with restart technique in parameter of regularization (see Chapter 3 [16]) allows us to reduce not strongly convex case to the strongly one (see also the text at once after the table in Part 1).

**Part 2. Projection Approach**

Consider the following problem with $g(x)$ of type 1 or 2 with $\lambda = 1$:

$$g(x) \to \min_{\substack{Ax=b \\ x \in Q}}. \qquad (3)$$

According to Lagrange principle we need solve the following saddle-point problem

$$\varphi(y) = \min_{x \in Q}\{g(x) + \langle y, Ax - b\rangle\} \to \max_y. \qquad (4)$$

We can solve (with $O(sn)$ flops) inner minimization problem explicitly and find

$$x(y), \ \varphi(y) = g(x(y)) + \langle y, Ax(y) - b\rangle, \ \nabla\varphi(y) = Ax(y) - b.$$

Function $g(x)$ is 1-strongly convex in $p$-norm, where $p = 2$ for the case 1 and $p = 1$ for the case 2. Hence $\varphi(y)$ has $L_2$-Lipschitz gradient in 2-norm ($L_2 = \sigma_{\max}(A)$ for the case 1 and $L_2 = \max_{k=1,\ldots,n} \|A^{\langle k \rangle}\|_2^2$ for the case 2 [12]). We use simple FGM [6], [12], [20] (not strongly convex variant) for dual problem (this dual problem (4) is different from the one in Remark 5) starting from $y^0 = 0$ of Nesterov's Universal method [22] (we hope that this method allows us to reduce $L_2$ to some extent) and with technique, described in items 2, 5.2,[6] 6.11 [15], we obtain (here we may use for FGM the estimation of accuracy certificate from [20], see formulas (3.3) – stopping criteria: this certificate is less then $\varepsilon$) after $O\left(sn\sqrt{L_2 R_{2,y}^2 / \varepsilon}\right)$ flops operations such $x$ that $\|Ax - b\|_2 \leq \varepsilon / R_{2,y}$, $g(x) - g_* \leq \varepsilon$, where $R_{2,y} = \|y_*\|_2$, $y_*$ – optimal solution of the dual problem (4), $x_*$ – optimal solution of the problem (3). This can be generalized if we use prox-FGM with $p \in [0,1]$ norm in $y$-space (with prox-function $d(x) = \frac{1}{p-1}\|x\|_p^2$ for $a \stackrel{def}{=} \frac{2\log m}{2\log m - 1} \leq p \leq 2$ and $d(x) = \frac{1}{a-1}\|x\|_a^2$ for $1 \leq p \leq a$): we obtain after $\tilde{O}\left(sn\sqrt{L_p R_{p,y}^2/\varepsilon}\right)$ flops operations[7] such $x$ that $\|Ax-b\|_q \leq \varepsilon/R_{p,y}$, $g(x) - g_* \leq \varepsilon$, where $1/p + 1/q = 1$, $R_{p,y} = \|y_*\|_p$. This estimation seems to be good enough according to Part 1. But unfortunately in typical applications $R_{p,y}$ can be large enough. If $g(x)$ has a bounded total variation $\Delta$ on $Q$ using technique of the Chapter 3 [16] we can change $b$ such that (we restrict ourselves here only the case $p = 2$) $\|b - \tilde{b}\|_2 \leq \delta$ and with this $\tilde{b}$ one have the following estimation $R_{2,y} = O(\Delta/\delta)$. So this approach though is fast enough but we suppose it isn't too fast as it could be according to the estimation above because of the $R_{2,y}$. Nevertheless we choose exactly this approach (to say more precisely we use close approach described in [21]).[8] The numerical results for real data we lead in the Part 3.

---

[6] In this item formula (46) a) can be write in a more precise way: $\|[Ax-b]^+\|_q \leq \varepsilon_{\text{cert}} / L$ – here we use denotations of [15]. In our paper this $L$ denoted by $R_{p,y}$.

**Remark.** These results (we restrict ourselves in this remark by the $p = 2$) can be also obtained based only on the Nesterov's estimated functions sequence technique [12]. The main ingredient of this approach (this can be generalized to composite stochastic optimization problem with inexact oracle Lemma 7.7 [16]) is such an estimation of convex functional $F(y)$ minimized on a convex set $\tilde{Q}$:

$$A_k F(\tilde{y}_k) \leq \min_{y \in \tilde{Q}} \Psi(y), \quad \Psi(y) = \beta_k d(y) + \sum_{i=0}^{k} a_i \{F(y_i) + \langle \nabla F(y_i), y - y_i \rangle\}, \quad A_k = \sum_{i=0}^{k} a_i.$$

From this one can easily find (it is significant here that $\tilde{Q}$ is a compact set, if this is not the truth one can consider instead of $\tilde{Q}$ for example the set $\tilde{Q} \cap \{y : d(y) \leq 2d(y_*)\}$, where $y_*$ – is a solution of the optimization problem)

$$F(\tilde{y}_k) - \min_{y \in \tilde{Q}} l(y) \leq \frac{\beta_k \max_{y \in \tilde{Q}} d(y)}{A_k} \stackrel{def}{=} \gamma_k, \quad l(y) = \sum_{i=0}^{k} a_i \{F(y_i) + \langle \nabla F(y_i), y - y_i \rangle\}.$$

For classical FGM $\beta_k = O(L_2)$, $a_k \sim k$, hence, $A_k \sim k^2$ and we have the rate of convergence $\sim k^{-2}$. Then one can use the technique of item 3 [23] to show that if $F(y) = \max_{x \in Q}\{\langle y, b - Ax \rangle - g(x)\}$ (denote $x(y)$ the solution of this problem) and the problem is $F(y) \to \min_{y \in \tilde{Q}}$ then

---

[7] Strictly speaking we have to use in this formula $R_{p,y}^2 \ln m$ instead of $R_{p,y}^2$ and with such prox-functions we have more expensive prox-mappings (at each iteration), because of impossibility of exact calculation according to explicit formulas. But it can be done no more than $O(sn \ln m) = \tilde{O}(sn)$ flops [14] (according to explicit formulas we need only $O(sn)$ flops). We should add that in this context we have to work with FGM with inexact oracle but not in conception of [16] where inexactness appears in gradient calculations. Our inexactness comes from prox-mapping, roughly speaking, when we calculate projection of the gradient. But it isn't very hard to generalize FGM on this conception of inexactness with example 4.1.3 a) [16] (as far as we know, this isn't done yet in general situation, but in our context one can find a precise calculation in item 5.5.1 [14] and [24] – it would be rather interesting to show how the results of these works can be obtained in conception of inexact oracle from Chapter 4 [16] but with inexactness depends (decrease) from the number of iteration $k$ as $\sim k^{-2}$).

[8] One can improve this approach. For this sake one has to consider randomized accelerated coordinate descent method for dual problem instead of simple FGM. This method, as we've already mentioned in Part 1, allows us to use sparcity of the matrix $A$ in a larger extent. This approach replaces worth-case Lipschitz constants to the average ones. Moreover this approach is also nice parallelized. But some technical moments arise in justification of its primal-dual nature. Nonetheless this can be done due to the work [20], since randomized accelerated coordinate descent method is a coupling of randomized non-accelerated coordinate descent method and primal-dual randomized coordinate descent mirror descent method (for smooth problems). So the desired properties are inherited from the mirror descent. Moreover if one will use regularization of dual problem [21] (to make it strongly convex) then the rate of convergence is linear and we have a convergence in argument, so there is no need in this regularized approach to use primal-dual nature of coordinate-descent.

$$0 \leq F(\tilde{y}_k) - \Phi(x_k) \leq \gamma_k, \quad \Phi(x) = -g(x) + \min_{y \in \tilde{Q}} \langle y, b - Ax \rangle, \quad x_k = \frac{1}{A_k} \sum_{i=0}^{k} a_i x(y_i).$$

If we use this result and choose $\tilde{Q} = \{y : \|y\|_2 \leq 2\tilde{R}_{2,y}\}$ (where $\tilde{R}_{2,y} = \max\left\{\sup_{k=0,1,2,\ldots} \|\tilde{y}_k\|_2, \|y_*\|_2\right\} = O(R_{2,y})$) we obtained

$$F(\tilde{y}_k) + g(x_k) + 2\tilde{R}_{2,y} \|Ax_k - b\|_2 \leq \gamma_k = O(k^{-2}).$$

From this result due to $Ax_* = b$ and the weak duality theorem $-g(x_*) = \Phi(x_*) \leq F(y_*)$ one can obtained

$$g(x_k) - g(x_*) \leq g(x_k) + F(y_*) \leq g(x_k) + F(\tilde{y}_k) \leq g(x_k) + F(\tilde{y}_k) + 2\tilde{R}_{2,y} \|Ax_k - b\|_2 \leq \gamma_k.$$

According to the definition of $F(y)$ and since $\|\tilde{y}_k\|_2 \leq \tilde{R}_{2,y}$, $\|y_*\|_2 \leq \tilde{R}_{2,y}$

$$-g(x_*) = \langle y_*, b - Ax_* \rangle - g(x_*) = F(y_*) \geq \langle y_*, b - Ax_k \rangle - g(x_k) \Rightarrow g(x_*) - g(x_k) \leq \tilde{R}_{2,y} \|Ax_k - b\|_2,$$

$$\tilde{R}_{2,y} \|Ax_k - b\|_2 \leq -g(x_k) + \langle \tilde{y}_k, b - Ax_k \rangle + g(x_k) + 2\tilde{R}_{2,y} \|Ax_k - b\|_2 \leq F(\tilde{y}_k) + g(x_k) + 2\tilde{R}_{2,y} \|Ax_k - b\|_2 \leq \gamma_k.$$

From the above one can conclude

$$|g(x_k) - g(x_*)| \leq \gamma_k, \quad \tilde{R}_{2,y} \|Ax_k - b\|_2 \leq \gamma_k.$$

One can mentioned that

$$g(x_k) - g(x_*) \leq F(\tilde{y}_k) + g(x_k) \leq \gamma_k.$$

So we have the following stopping rule condition (we choose $\varepsilon$ and $\tilde{\varepsilon}$ in advance)

$$F(\tilde{y}_k) + g(x_k) \leq \varepsilon \text{ and } \|Ax_k - b\|_2 \leq \tilde{\varepsilon}. \qquad (*)$$

From the above we can guarantee that (*) is fulfill for

$$k = \max\left\{\sqrt{\frac{8L_2 \tilde{R}_{2,y}^2}{\varepsilon}}, \sqrt{\frac{8L_2 \tilde{R}_{2,y}}{\tilde{\varepsilon}}}\right\}.$$

We can consider unconstrained variant of the dual problem (that is the real dual problem): $F(y) \to \min_y$. Since $\tilde{R}_{2,y}$ isn't included in stepsize policy of the FGM for $F(y) \to \min_{y \in \tilde{Q}}$ we have that both of the problems generate the same sequences if $\tilde{R}_{2,y}$ is taken sufficiently large. Therefore we can restrict ourselves only by consideration of real dual problem with FGM (that is we use FGM with $\tilde{Q} = \mathbb{R}^m$) and stopping rule condition (*). If we want to estimate $\tilde{R}_{2,y}$ (for example we want to choose an optimal correspondence between $\varepsilon$ and $\tilde{\varepsilon}$) we can use here the restart technique: we fixed $\tilde{R}_{2,y}$ and make, according to this $\tilde{R}_{2,y}$, prescribed number of iterations (we have explicit formula for that). Then we can verify $\|\tilde{y}_k\|_2 \leq \tilde{R}_{2,y}$ and stopping rule conditions (*). If they are fulfilled we stop. If this is not the true, we put $\tilde{R}_{2,y} := 2\tilde{R}_{2,y}$ and restart procedure. One can show that restart constant 2 here is optimal (see the end of item 6 [21]).

All the results can be generalized for the case of inequalities $Ax \leq b$ in (3) instead of $Ax = b$ (see for example item 5.2 [15]). Moreover, we can consider more general 1-strongly convex functions $g(x)$ such that inner problem in (4) hasn't explicit solution. In this case we can use item 4.2.2 of [16] (see example 4 in [25]). If we don't have strong convexity of $g(x)$ we may use separable structure of $g(x)$. In case when $Q$ is also separable, i.e. parallelepiped (or can consists of a few simple constraints, for example simplex type, Euclidian ball type (see text after Remark 4 in Part 1)), we can decompose prox-mapping problem to $n$ one-dimensional convex problems. Each of these problems can be efficiently solve with linear rate of convergence what is the clue to success. After that we may use some subgradient type method for outer (not smooth in general) convex problem with conception of $\varepsilon$-subgradient (see item 5 of paragraph 1 Chapter 5 [5]).

Note that according to the penalty functions method (see Theorem 7 in section 2.5 Chapter 8 [5] and section 16 Chapter 5 [26]) we have (assume that $A$ is a full-rank matrix) that for the solution $x_K$ of the problem

$$g(x) + \frac{K}{2}\|Ax - b\|_2^2 \to \min_{x \in Q},$$

for $K \to \infty$: $\|x_K - x_*\| = O(1/K)$, $K \cdot (Ax_K - b) \to y_*$, where $x_*$ is solution of the primal problem (3) and $y_*$ is solution of the dual problem (4). This fact allows us to connect the problems formulations in part 1 and part 2 (without any tricks described below).

Now let's return to the problem of choosing of structural parameter mentioned in Part 1. We consider here cases 1 and 2 (case 3 see in [1]). The idea is as follows (we demonstrate the idea on the case 1). Consider the 1-dimensional smooth convex optimization problem (see Introduction)

$$\min_{x \geq 0}\left\{\|x - x_g\|_2^2 + \frac{\bar{\lambda}}{2}\left(\|Ax - b\|_2^2 - \varepsilon^2\right)\right\} \to \max_{\bar{\lambda} \geq 0}.$$

If we find the solution of this problem we automatically find structural parameter $\lambda = \bar{\lambda}^{-1}$. Since we can solve with fixed $\bar{\lambda}$ inner problem with linear rate of convergence (see Part 1), we can calculate approximate $\lambda$-derivative and use concept of inexact oracle from [16] or $\varepsilon$-subgradient (see above) to solve the outer problem on $\lambda$. But instead of FGM (see example 4 in [25]) we use golden section search for outer problem (see, for example, [5]) with linear rate of convergence independently of any strong convexity. So it can be shown that problems (1) and (3) are $\varepsilon$-equivalent up to logarithmic factor.

Let's consider at the end one interesting special example. Assume that we know in advance that the 1-norm of the optimal dual solution in (4) is $R_{1,y}$. In this case we can reformulate the problem (4) as

$$\|Ax - b\|_\infty R_{1,y} + g(x) \to \min_{x \in Q},$$

$$\|Ax - b\|_\infty + \frac{1}{R_{1,y}} g(x) \to \min_{x \in Q}.$$

And similarly

$$\|Ax - b\|_q + \frac{1}{R_{p,y}} g(x) \to \min_{x \in Q}, \quad \frac{1}{q} + \frac{1}{p} = 1.$$

These problems aren't smooth and, as a consequence, aren't such computational simple as, for example, problem (1). Exclusion is the case 3 $g(x) = \|x\|_1$, $p = 1$, $Q = \mathbb{R}_+^n$. In this case we can reformulate one more time this problem to reduce it

$$\|Ax - b\|_\infty \to \min_{x \in \bar{S}_n(R(R_{1,y}))}.$$

The last problem can be solved by randomized mirror descent [27] (see item 6.5.2.3). The only problem is $R(R_{1,y})$ – we don't know this constant as a rule.

**Part 3. Results of numerical experiment with real data**

We provide experiments on simulated data. In all the cases we test only fast gradient algorithm, described in Part 2, being applied to the dual problem (for Mimimal mutual information model).

In the first experiment we have a flat network topology with 1000 nodes, $10^6$ demands and 10000 links. We suppose that demand value between each pair of nodes is a uniformly distributed in [100, 300] random variable.

We refer link load estimation quality as $LLA = \|Ax - b\|_2 / \|b\|_2$ and demands estimation quality as $DA = \|x - x'\|_2 / \|x'\|_2$, where $x'$ is a true solution. Solution time is 2.132 seconds on 3.7 Gz machine with 8 Gb RAM, demand estimation quality $DA$ is 66.33%.

In second experiment we have a flat network topology with 2000 nodes, $4 \cdot 10^6$ demands and 20000 links. We suppose that demand value between each pair of nodes is a uniformly distributed in [100, 300] random variable. For this topology $DA = 66.36\%$ and $LLA = 99\%$. Solution time is 7.345 seconds on the same machine.

We'd like to thanks prof. A. Nemirovski, prof. Yu. Nesterov and prof. B. Polyak for fruitful discussions.

The work was partially supported by RFBR (projects 14-01-00722-a, 15-31-20571-mol-a-ved) and partially supported by Russian Scientific Fund (project 14-50-00150).

**References**

[1]. *Hastie T., Tibshirani R., Friedman R.* The Elements of statistical learning: Data mining, Inference and Prediction. Springer, 2009. http://statweb.stanford.edu/~tibs/ElemStatLearn/

[2]. *Spokoiny V.G.* Penalized maximum likelihood estimation and effective dimension // e-print, 2012. arXiv:1205.0498

[3]. *Zhang Y., Roughan M., Duffield N., Greenberg A.* Fast Accurate Computation of Large-Scale IP Traffic Matrices from Link Loads // In ACM Sigmetrics. San Diego, CA, 2003. https://www.cs.utexas.edu/~yzhang/papers/tomogravity-sigm03.pdf

[4]. *Zhang Y., Roughan M., Lund C., Donoho D.* Estimating Point-to-Point and Point-to-Multipoint Traffic Matrices: An Information-Theoretic Approach // IEEE/ACM Transactions of Networking. 2004. V. 10. № 10. https://www.cs.utexas.edu/~yzhang/papers/mmi-ton05.pdf

[5]. *Polyak B.T.* Introduction to optimization. Hardcover, 1987.

[6]. *Nesterov Yu.* Gradient methods for minimizing composite functions // Math. Prog. 2013. V. 140. № 1. P. 125–161.

[7]. *Fercoq O., Richtarik P.* Accelerated, Parallel and Proximal Coordinate Descent // e-print, 2013. arXiv:1312.5799

[8]. *Qu Z., Richtarik P.* Coordinate Descent with Arbitrary Sampling I: Algorithms and Complexity // e-print, 2014. arXiv:1412.8060

[9]. *Shalev-Shwartz S., Zhang T.* Accelerated proximal stochastic dual coordinate ascent for regularized loss minimization // In Proceedings of the 31th International Conference on Machine Learning, ICML 2014, Beijing, China, 21-26 June 2014. P. 64–72. arXiv:1309.2375

[10]. *Zheng Q., Richtárik P., Zhang T.* Randomized dual coordinate ascent with arbitrary sampling // e-print, 2015. arXiv:1411.5873

[11]. *Beck A., Teboulle M.* A Fast Iterative Thresholding Algorithm for Linear Inverse Problem // SIAM Journal of Image Sciences. 2009. V. 2(1). P. 183–202.

[12]. *Nesterov Y.* Smooth minimization of non-smooth function // Math. Program. Ser. A. 2005. V. 103. № 1. P. 127–152.

[13]. *Juditsky A., Nesterov Yu.* Deterministic and stochastic primal-dual subgradient algorithms for uniformly convex minimization // Stoch. System. 2014. V. 4. no. 1. P. 44–80. arXiv:1401.1792

[14]. *Nemirovski A.* Lectures on modern convex optimization analysis, algorithms, and engineering applications. Philadelphia: SIAM, 2013. http://www2.isye.gatech.edu/~nemirovs/Lect_ModConvOpt.pdf

[15]. *Nemirovski A., Onn S., Rothblum U.* Accuracy certificates for computational problems with convex structure // Mathematics of Operations Research. 2010. V. 35:1. P. 52–78.

[16]. *Devolder O.* Exactness, inexactness and stochasticity in first-order methods for large-scale convex optimization. CORE UCL, PhD thesis, March 2013.

[17]. *Gasnikov A., Dvurechensky P., Usmanova I.* About accelerated randomized methods // TRUDY MIPT. 2016. V. 8. no. 2. (in print) arXiv:1508.02182 [in Russian]

[18]. *Boyd S., Parikh N., Chu E., Peleato B., Eckstein J.* Distributed optimization and statistical learning via the alternating direction method of multipliers // Foundations and Trends in Machine Learning. 2011. V. 3(1). P. 1–122. http://stanford.edu/~boyd/papers.html

[19]. *Richtárik P.* http://www.maths.ed.ac.uk/~richtarik/

[20]. *Allen-Zhu Z., Orecchia L.* Linear coupling: An ultimate unification of gradient and mirror descent // e-print, 2014. arXiv:1407.1537

[21]. *Gasikov A.V., Gasnikova E.V., Nesterov Yu.E., Chernov A.V.* Entropy-linear programming // Comp. Math. and Math. Phys. 2016. V. 56. no. 4. P. 523–534. arXiv:1410.7719

[22]. *Nesterov Yu.* Universal gradient methods for convex optimization problems // CORE Discussion Paper 2013/63. 2013.

[23]. *Nesterov Yu.* Complexity bounds or primal-dual methods minimizing the model of objective function // CORE Discussion Paper 2015/3.
http://www.uclouvain.be/cps/ucl/doc/core/documents/coredp2015_3web.pdf

[24]. *Nesterov Yu., Nemirovski A.* On first order algorithms for $l_1$ / nuclear norm minimization // Acta Numerica. 2013. V. 22. P. 509–575.

[25]. *Gasnikov A., Dvurechensky P., Nesterov Yu.* Stochastic gradient methods with inexact oracle // TRUDY MIPT. 2016. V. 8. no. 1. P. 41–91. arxiv:1411.4218 [in Russian]

[26]. *Vasiliev F.P.* Optimization methods. M. MCCME, 2011. Vol. 1. [in Russian]

[27]. *Juditsky A., Nemirovski A.* First order methods for nonsmooth convex large-scale optimization, I, II. In: Optimization for Machine Learning. Eds. S. Sra, S. Nowozin, S. Wright. MIT Press, 2012.